\documentclass[journal]{IEEEtran}
\usepackage{amssymb}
\usepackage[pdftex]{graphicx}
\usepackage{tikz}

\usepackage{framed}
\usepackage{bm}
\usepackage{cite}
\usepackage{comment}
\usepackage{enumerate}
\usepackage{enumitem}
\usepackage{amsmath}
\usepackage{amsfonts}
\usepackage{xcolor}
\usepackage{cases}
\usepackage{mathtools}
\mathtoolsset{showonlyrefs=true}

\usetikzlibrary{intersections,calc,arrows.meta}

\newtheorem{theorem}{Theorem}

\newtheorem{remark}{Remark}
\newtheorem{lemma}{Lemma}
\newtheorem{corollary}{Corollary}

\newtheorem{proposition}{Proposition}

\newtheorem{example}{Example}
\usepackage{algpseudocode}
\newcounter{algorithm}
\newenvironment{algorithm}[1][]{%
  \par\addvspace{0.8\baselineskip}\noindent\begin{minipage}{\columnwidth}%
  \refstepcounter{algorithm}\hrule height 0.8pt\relax
  \vspace{0.4em}\noindent\textbf{Algorithm \thealgorithm. #1}\par
  \vspace{0.4em}\hrule height 0.4pt\relax\vspace{0.2em}%
}{\vspace{0.2em}\hrule height 0.8pt\relax\end{minipage}\par\addvspace{0.8\baselineskip}}
\algrenewcommand{\algorithmicrequire}{\textbf{Input:}}
\algrenewcommand{\algorithmicensure}{\textbf{Output:}}

\usepackage{blkarray}

\usepackage[hidelinks]{hyperref}
\usepackage[sort]{cleveref}
\usepackage{mathtools}
\usepackage{multirow}
\usepackage{threeparttable}

\usepackage{pgfplots}
\pgfplotsset{compat=1.18}
\usepgfplotslibrary{statistics}

\crefname{section}{Section}{Sections}
\crefname{appendix}{Appendix}{Appendices}
\crefname{theorem}{Theorem}{Theorems}
\crefname{proposition}{Proposition}{Propositions}
\crefname{lemma}{Lemma}{Lemmas}
\crefname{problem}{Problem}{Problems}
\crefname{assumption}{Assumption}{Assumptions}
\crefname{algorithm}{Algorithm}{Algorithms}
\crefname{remark}{Remark}{Remarks}
\crefname{equation}{}{}
\Crefname{equation}{Equation}{Equations}
\crefname{example}{Example}{Examples}
\crefname{figure}{Fig.}{Figs.}
\crefname{table}{Table}{Tables}

\crefname{enumi}{}{}

\DeclareMathOperator*{\argmin}{arg\,min}

\DeclareMathOperator{\diag}{diag}

\DeclareMathOperator{\tr}{tr}

\newcommand{\R}{\mathbb R}
\newcommand{\one}{\bm 1}
\newcommand{\vecop}{\operatorname{vec}}
\newcommand{\CS}{\mathrm{CS}}

\newcommand{\sgn}{\operatorname{sgn}}

\usepackage{latexsym}

\usepackage{url}
\definecolor{RankFirst}{HTML}{005A9C}
\definecolor{RankSecond}{HTML}{9B4800}
\definecolor{RankThird}{HTML}{713B86}

\begin{document}
\title{Robustness of controllability scores\\ for stable linear network systems}
\author{Ritsuki~Yamada and Kazuhiro~Sato\thanks{R. Yamada and K. Sato are with the Department of Mathematical Informatics, Graduate School of Information Science and Technology, The University of Tokyo, Tokyo 113-8656, Japan.}}

\maketitle
\thispagestyle{empty}
\pagestyle{empty}
\lefthyphenmin=2\righthyphenmin=3
\begin{abstract}
Controllability scores quantify node importance in
linear network systems by optimizing the allocation
of input strength. Although the volumetric
controllability score (VCS) and average-energy
controllability score (AECS) are uniquely defined
for Hurwitz systems, their sensitivity to
perturbations of the system matrix remains unclear.
In this paper, we investigate the robustness
of these scores to such perturbations.
We first prove that both scores are locally
Lipschitz continuous in the system matrix and establish
a robustness bound for Hurwitz comparison matrices.
We then characterize the directional derivative
of the optimal allocation along a differentiable
matrix path as the unique solution of a
critical-cone quadratic program, whose linear
term is computed using differentiated Lyapunov
equations.
For single-edge weakening, the norm of the
resulting response vector provides a
robustness-based measure of edge importance.
Numerical experiments show that the response, and hence the
edge importance, depends on both the criterion and the dynamics, not
only on the network topology.
\end{abstract}
\IEEEpeerreviewmaketitle
\section{Introduction}
\label{sec:introduction}
Centrality measures provide a way to compare the importance of nodes in a network. Classical measures distinguish different structural roles, such as direct connectivity, proximity to other nodes, and mediation of paths~\cite{freeman1979centrality,rodrigues2019centrality}, while different flow mechanisms can support different centrality measures~\cite{borgatti2005centrality}. For network control, this motivates centralities that account for the dynamics and the chosen control objective, rather than on topology alone~\cite{pasqualetti2014metrics,summers2016submodularity,sato2024scores}.

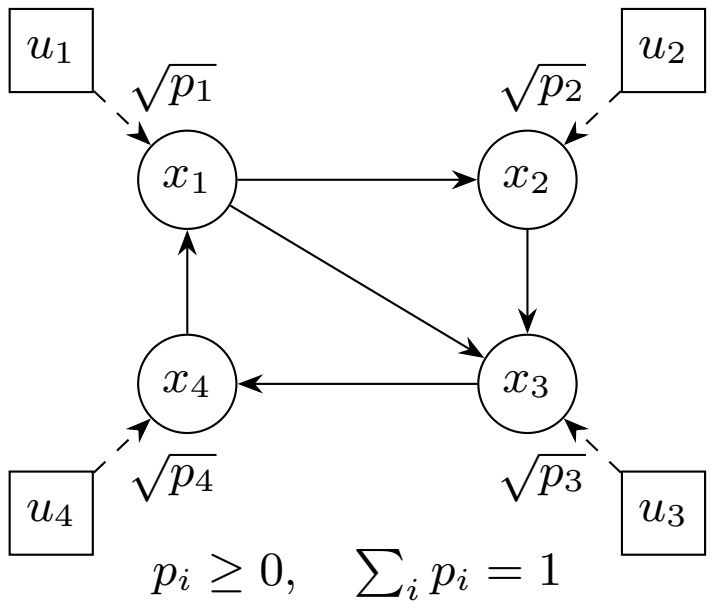
\begin{figure}[t]
\centering
\begin{tikzpicture}[>=Stealth,scale=0.9,every node/.style={font=\small},
 state/.style={draw,circle,minimum size=6.5mm,inner sep=1pt},
 input/.style={draw,rectangle,minimum size=5.5mm,inner sep=1pt}]
\node[state] (x1) at (0,0) {$x_1$};
\node[state] (x2) at (2.5,0) {$x_2$};
\node[state] (x3) at (2.5,-1.5) {$x_3$};
\node[state] (x4) at (0,-1.5) {$x_4$};
\node[input] (u1) at (-1,0.95) {$u_1$};
\node[input] (u2) at (3.5,0.95) {$u_2$};
\node[input] (u3) at (3.5,-2.45) {$u_3$};
\node[input] (u4) at (-1,-2.45) {$u_4$};
\draw[->] (x1) -- (x2);
\draw[->] (x2) -- (x3);
\draw[->] (x3) -- (x4);
\draw[->] (x4) -- (x1);
\draw[->] (x1) -- (x3);
\draw[->,dashed] (u1) -- node[above right,inner sep=1pt] {$\sqrt{p_1}$} (x1);
\draw[->,dashed] (u2) -- node[above left,inner sep=1pt] {$\sqrt{p_2}$} (x2);
\draw[->,dashed] (u3) -- node[below left,inner sep=1pt] {$\sqrt{p_3}$} (x3);
\draw[->,dashed] (u4) -- node[below right,inner sep=1pt] {$\sqrt{p_4}$} (x4);
\node at (1.25,-2.9) {$p_i\geq0,\quad\sum_i p_i=1$};
\end{tikzpicture}
\caption{Virtual node-wise inputs used to define VCS and AECS. Circles represent states, squares represent inputs, solid arrows represent network couplings, and dashed arrows represent input gains $\sqrt{p_i}$. The allocation $p$ is optimized under $p_i\geq0$ and $\sum_i p_i=1$.}
\label{fig:score-construction}
\end{figure}

Controllability scores~\cite{sato2024scores} are one such centrality. For autonomous network dynamics described by
\begin{equation}\label{eq:autonomous}
 \dot x=Ax,
\end{equation}
they compare hypothetical inputs applied separately at all nodes, as illustrated in Fig.~\ref{fig:score-construction}. The system matrix $A$ encodes the weighted network structure. A nonnegative variable $p_i$ specifies the squared input strength at node $i$, with $\sum_i p_i=1$. Optimizing this allocation defines node centralities. The volumetric controllability score (VCS) maximizes the volume of the ellipsoid associated with the controllability Gramian, whereas the average-energy controllability score (AECS) minimizes an average minimum-energy criterion~\cite{sato2024scores}. When~\eqref{eq:autonomous} is asymptotically stable, both infinite-horizon optimal allocations exist and are unique. Each criterion therefore assigns one well-defined vector of node centralities to network system~\eqref{eq:autonomous}.

Uniqueness makes the scores unambiguous at a fixed model, but it does not establish how they change when that model is perturbed. Connection weights estimated from data may be inaccurate, and the connections themselves may vary~\cite{arola2020uncertainty,holme2012temporal}. To inform the choice of control nodes, the scores should therefore be interpreted together with their sensitivity to these changes. Segarra and Ribeiro formalized this issue for weighted-graph centrality measures and showed that robustness distinguishes their responses to noisy weights~\cite{segarra2016stability}. Because VCS and AECS also depend on the dynamics and control objective~\cite{sato2024scores}, their robustness, i.e., score changes under model perturbations, must be examined through the corresponding scoring problem.

We address this question at two levels, bounds on score changes and the first-order response to a specified perturbation, and use the latter to assess edge importance. The contributions of this paper are summarized as follows.
\begin{enumerate}[leftmargin=*,label=\arabic*)]
\item We prove that VCS and AECS are locally Lipschitz continuous in the system matrix $A$. Combined with the boundedness of the simplex, this yields a reference-dependent bound on score changes, analogous to centrality robustness for weighted graphs~\cite{segarra2016stability}.

\item We characterize the one-sided derivative of the scores along a
prescribed differentiable matrix path as the unique solution of a
critical-cone quadratic program. The resulting controllability score
response vector describes which node scores increase or decrease, to
first order, under a specified model change. It is computable from
differentiated Lyapunov equations and one convex QP per direction,
without reoptimizing the allocation at nearby models.

\item For single-edge weakening, the norm of the response vector
defines a robustness-based edge importance. Edges with a large norm
can be interpreted as those to which the scores are most sensitive,
so the accuracy of their estimated weights matters for the
reliability of the scores.

\item We numerically illustrate that the response depends on the dynamics, not only on
the network topology. For the same graph and the same weakened edge,
Laplacian and fixed-diagonal dynamics give different response norms,
and some node-wise components have opposite signs. Thus, edge
importance for controllability scores cannot be inferred from the
network topology alone, as in structural approaches~\cite{liu2011controllability,mengiste2015edge,chung2021dilations}, so it is useful
to report the dynamical model together with the edge importance.
\end{enumerate}

The rest of this paper is organized as follows.
Section~\ref{sec:formulation} formulates the scores and perturbation
paths. Section~\ref{sec:continuity} gives the robustness bounds of
contribution~1). Section~\ref{sec:sensitivity} derives the response
vector of contribution~2) and defines the edge importance of contribution~3). Section~\ref{sec:numerical} examines the dependence of the responses on the criterion and dynamics, as stated in contribution~4), and validates the accuracy of the first-order predictions. Appendices~\ref{app:lipschitz} and~\ref{app:derivative} contain the proofs, Appendix~\ref{app:comparison} gives the analytical comparison of the two criteria and Appendix~\ref{app:metzler}
gives a full-path admissibility result for Metzler systems.

\emph{Notation:} The sets of real numbers, real $n$-vectors, and real
$n\times n$ matrices are $\mathbb R$, $\mathbb R^n$, and
$\mathbb R^{n\times n}$, and $\mathbb S^n$ is the space of real
symmetric $n\times n$ matrices. All vectors are column vectors. The
vector $\one$ has all components equal to one, $I$ is the identity
matrix, $e_i$ is the $i$th coordinate vector, $E_{rs}=e_re_s^\top$, and $(\cdot)^\top$ denotes transpose. For a vector $v$,
$v^\perp=\{z:v^\top z=0\}$. For an index set $J$, $J^c$ is its
complement, $z_J$ is the subvector of $z$ indexed by $J$, and $H_{JJ}$
is the corresponding principal submatrix of $H$. We write $\tr$ and
$\det$ for the trace and determinant, $\lambda_{\min}$ and
$\lambda_{\max}$ for the smallest and largest eigenvalues of a
symmetric matrix, $\sigma_{\min}$ for the smallest singular value,
$\diag$ for diagonal-matrix formation, $\vecop$ for vectorization, and
$\langle U,V\rangle_{\mathrm F}=\tr(U^\top V)$. The norm
$\lVert\cdot\rVert_2$ is the Euclidean norm for vectors and the
spectral norm for matrices, and $\lVert\cdot\rVert_{\mathrm F}$ is the
Frobenius norm. The relations $\succ$ and $\succeq$ denote the
positive-definite and positive-semidefinite orders; entrywise
inequalities are stated explicitly. For a function $f(p,A)$,
$\nabla_p f$ and $\nabla^2_{pp}f$ are its gradient and Hessian in $p$.
As $\theta\downarrow0$, $o(\theta)$ and $O(\theta)$ are the usual
Landau symbols, and $\sgn$ is the sign function.

\section{Preliminaries and problem formulation}
This section reviews the infinite-horizon VCS and AECS and their optimization formulations, and then introduces the graph distance and matrix perturbation paths used to analyze how these scores change under model perturbations.

Throughout this paper, we assume that the system matrix $A\in\mathbb R^{n\times n}$, $n\geq2$, in~\eqref{eq:autonomous} is Hurwitz, that is, all its eigenvalues have negative real parts.

\label{sec:formulation}
\subsection{Infinite-horizon controllability scores}
To define its controllability scores for system~\eqref{eq:autonomous}, we augment the system with virtual node-wise inputs:
\begin{equation}\label{eq:system}
 \dot x=Ax+B(p)u,\qquad
 B(p)=\diag(\sqrt{p_1},\ldots,\sqrt{p_n}),
\end{equation}
where $p$ belongs to $\Delta=\{p\in\mathbb R^n:p\geq0,\one^\top p=1\}$. Thus, $p_i$ is the share of the squared input strength assigned to node $i$, and $p$ is a continuous allocation.

In order to measure the degree of controllability, we use the controllability Gramian~$W_A(p)$ defined as
\begin{equation}\label{eq:WA}
 W_A(p)=\sum_{i=1}^n p_iW_i(A),
\end{equation}
where
\begin{equation}\label{eq:Wi}
 W_i(A)=\int_0^\infty e^{At}e_ie_i^\top e^{A^\top t}\,dt.
\end{equation}
Because $A$ is Hurwitz, $e^{At}$ decays exponentially, so the integrals in~\eqref{eq:Wi} converge. In particular, each $W_i(A)$ satisfies the Lyapunov equation
\begin{equation}\label{eq:lyapunov}
 AW_i(A)+W_i(A)A^\top+e_ie_i^\top=0.
\end{equation}
The admissible set is $X_A=\{p\in\Delta:W_A(p)\succ0\}$. The uniform allocation $p^0=\one/n$ belongs to $X_A$, since $B(p_0)=I/\sqrt{n}$ is nonsingular and hence $W_A(p_0) \succ 0$.

The two scores optimize different scalar summaries of this Gramian. Following~\cite{sato2024scores}, for $\CS\in\{\mathrm V,\mathrm A\}$, define
\begin{equation}\label{eq:score}
 p^{\CS}(A)=\argmin_{p\in X_A}f_{\CS}(p,A),
\end{equation}
\begin{equation}\label{eq:objectives}
 \begin{aligned}
 f_{\mathrm V}(p,A)&=-\log\det W_A(p),\\
 f_{\mathrm A}(p,A)&=\tr\bigl(W_A(p)^{-1}\bigr).
 \end{aligned}
\end{equation}
The minimizing vectors $p^{\mathrm V}(A)$ and $p^{\mathrm A}(A)$ are the VCS and AECS. For the Hurwitz, node-wise-input formulation above, both minimizing allocations exist and are unique~\cite[Sec.~III]{sato2024scores}. This makes each component a well-defined measure of relative node importance: it is the share of the squared input-strength budget assigned to that node under the chosen criterion.

For a positive-definite Gramian $W_A(p)$, the ellipsoid $\{x\in\mathbb R^n:x^\top W_A(p)^{-1}x\leq1\}$ has volume proportional to $\det(W_A(p))^{1/2}$, and the mean quadratic energy over uniformly distributed unit target directions is $\tr(W_A(p)^{-1})/n$ \cite{muller1972analysis}. These standard Gramian interpretations motivate maximizing the volume and minimizing the average minimum energy.

The perturbation analysis uses the gradient and Hessian of $f_{\CS}$
with respect to the allocation $p$. The functions $f_{\mathrm V}$ and $f_{\mathrm A}$ are smooth in $p$ wherever $W_A(p)\succ0$. The gradient and Hessian formulas for VCS
and AECS are given in~\cite[Sec.~III]{sato2024scores}.
At a feasible $p$, let $R=W_A(p)^{-1}$ and $H_{\CS}=\nabla^2_{pp}f_{\CS}$. The gradient components are
\begin{equation}\label{eq:gradients}
 (\nabla_p f_{\mathrm V})_i=-\tr(RW_i),\qquad
 (\nabla_p f_{\mathrm A})_i=-\tr(RW_iR).
\end{equation}
For
\begin{equation}\label{eq:QA}
Q_A(z)=\sum_i z_iW_i(A),
\end{equation}
the Hessian quadratic forms are
\begin{equation}\label{eq:HV}
 z^\top H_{\mathrm V}z=\tr\bigl(RQ_A(z)RQ_A(z)\bigr),
\end{equation}
\begin{equation}\label{eq:HA}
 z^\top H_{\mathrm A}z=2\tr\bigl(RQ_A(z)RQ_A(z)R\bigr).
\end{equation}

\begin{remark}[Other scoring formulations]
There are also guarantees beyond Hurwitz dynamics, under different formulations. At a finite horizon, the controllability-scoring problems have minimizers for every real $A$, and these minimizers are unique for almost every terminal time $T>0$~\cite[Thm.~1]{sato2025uniqueness}. For non-Hurwitz systems, \cite{umezu2026infinite} defines limiting controllability scores as the horizon tends to infinity using scaled Gramians. The limiting VCS accounts for all spectral components, whereas the limiting AECS depends only on the stable modes and requires at least one stable eigenvalue. Under suitable rank conditions, these limiting scores are unique. The present perturbation analysis concerns the conventional unscaled Hurwitz formulation in~\eqref{eq:score}--\eqref{eq:objectives}.
\end{remark}

\subsection{Graph distance and edge perturbations}
To compare the scores across models, we specify which matrix entries may change and how their changes are measured. Our convention is that $A_{rs}$ represents the directed connection $s\to r$. All comparisons use the same node labeling. Let $\mathcal E$ be a fixed set of allowed matrix positions. For \(A,\widetilde A\in\mathbb R^{n\times n}\) satisfying
$A_{(r,s)}=\widetilde A_{(r,s)}, \quad \text(r,s)\notin\mathcal E$,
define
\begin{equation}\label{eq:graph-distance}
 d_{\mathcal E}(A,\widetilde A)
   =\sum_{(r,s)\in\mathcal E}|A_{rs}-\widetilde A_{rs}|.
\end{equation}
Zero weights are allowed, and the distance therefore covers edge insertion and deletion. For such matrices, $\lVert A-\widetilde A\rVert_{\mathrm F}\leq d_{\mathcal E}(A,\widetilde A)$.

Whereas the distance measures the size of a model change, a matrix path also specifies its direction $D \in \R^{n\times n}$. We consider paths defined for sufficiently small $\theta\geq0$, with
\begin{equation}\label{eq:matrix-path}
 A(\theta)=A+\theta D+o(\theta).
\end{equation}
Since $A$ is Hurwitz, $A(\theta)$ remains Hurwitz near $\theta=0$, and $D$ can specify simultaneous changes at several matrix positions.

\begin{example}[Edge weakening]\label{ex:edge_weakening}
Let $A_{rs}\ne0$ and $r\ne s$. Weakening this entry while keeping the others fixed gives
\begin{equation}\label{eq:deletion}
 A(\theta)=A-\theta A_{rs}E_{rs},\qquad0\leq\theta\leq1,
\end{equation}
with $D=-A_{rs}E_{rs}$. Here $\theta=0$ is the reference and $\theta=1$ removes the entry. The absolute change is $\varepsilon=\theta|A_{rs}|=\lVert A(\theta)-A\rVert_{\mathrm F}$.
\end{example}

\begin{example}[Row-Laplacian edge weakening]\label{ex:laplacian}
Consider the row-Laplacian model $A=-L-\alpha I$ with $\alpha>0$,
where $\mathrm{Adj}\in\mathbb R^{n\times n}$ is a nonnegative weighted
adjacency matrix, with $\mathrm{Adj}_{sr}$ the weight of the edge
$s\to r$, and
\[
 L=\diag(\ell_1,\ldots,\ell_n)-\mathrm{Adj}^\top,\qquad
 \ell_r=\sum_s \mathrm{Adj}_{sr},
\]
so that $L\one=0$. This model adds a uniform leakage term $-\alpha x$
to the first-order consensus dynamics
$\dot x=-Lx$~\cite{olfati2007consensus}; the leakage makes $A$
Hurwitz, whereas $-L$ alone has a zero eigenvalue. For undirected
networks, $L+\alpha I$ is a loopy Laplacian, as arises in RC circuits
with uniform shunts to ground~\cite{dorfler2018electrical}. Weakening a positive edge weight $w=\operatorname{Adj}_{sr}$ in this model also
reduces $\ell_r$ by the same amount. The resulting path is
\[
 A(\theta)=A+\theta w(E_{rr}-E_{rs}),\qquad 0\leq\theta\leq1,
\]
so $D=w(E_{rr}-E_{rs})$ includes the associated diagonal change,
unlike the fixed-entry path of Example~1.
\end{example}

A useful subclass is given by Metzler matrices, whose off-diagonal entries are nonnegative. They describe continuous-time positive systems and arise, for example, in transportation networks and vehicle-formation models~\cite{rantzer2015scalable}, as well as in the row-Laplacian
model of Example~\ref{ex:laplacian}. For a Hurwitz Metzler matrix, the fixed-entry deletion path~\eqref{eq:deletion} remains Hurwitz for all \(0\le\theta\le1\), as proved in Appendix~\ref{app:metzler}.

\section{Robustness under model perturbations}
\label{sec:continuity}
In this section, we bound how much the VCS and AECS can change
when the system matrix is perturbed. The argument has two steps.
Lemma~\ref{lem:curvature} gives a lower bound on the Hessian quadratic form of each objective~$f_{\rm CS}$
along directions $z \in \one^\perp$ uniformly
over $p\in X_A$. Lemma~\ref{lem:lipschitz} combines this curvature with a
perturbation bound on the gradient to show that the optimal
allocation~$p^{\rm CS}$ is locally Lipschitz continuous in the system matrix~$A$.
Theorem~\ref{thm:graph} then extends this local estimate to a bound relative to a
fixed reference matrix.

The first step is a curvature bound that does not depend on the
allocation. It uses the injectivity of the Gramian map established
in~\cite{sato2024scores} and has the same quantitative form as the finite-horizon
restricted-curvature bound in~\cite[Lem.~6]{sato2026projected}.

\begin{lemma}[Hessian lower bound on $\one^\perp$]\label{lem:curvature}
For every Hurwitz matrix $A$ and criterion $\CS$, there is a constant $m_{\CS}(A)>0$, independent of $p\in X_A$, such that
\begin{equation}\label{eq:curvature}
 z^\top H_{\CS}(p,A)z\geq m_{\CS}(A)\lVert z\rVert_2^2,
 \qquad z\in\one^\perp.
\end{equation}
\end{lemma}
\begin{IEEEproof}
See Appendix~\ref{app:curvature}.
\end{IEEEproof}

We next establish the local Lipschitz continuity of the optimal allocation. Since the admissible set \(X_A\) depends on \(A\), comparing the optimality conditions at two nearby matrices requires showing that each optimizer is feasible at the other matrix. Combining this property with the curvature bound of Lemma~\ref{lem:curvature} and a gradient perturbation bound yields the following result.

\begin{lemma}[Local Lipschitz continuity]\label{lem:lipschitz}
Fix a Hurwitz matrix \(A_0\) and a criterion \(\mathrm{CS}\in\{\mathrm V,\mathrm A\}\). There exists an open neighborhood $U$ of $A_0$, consisting of Hurwitz matrices, and $L_{\CS}>0$ such that
\begin{equation}\label{eq:lipschitz}
 \lVert p^{\CS}(A)-p^{\CS}(\widetilde A)\rVert_2
 \leq L_{\CS}\lVert A-\widetilde A\rVert_{\mathrm F}\quad (A,\widetilde A\in U).
\end{equation}
\end{lemma}
\begin{IEEEproof}
See Appendix~\ref{app:lipschitz}.
\end{IEEEproof}

Lemma~\ref{lem:lipschitz} applies also to boundary allocations, without assuming strict complementarity or that the support of the optimal allocation remains fixed.

Lemma~\ref{lem:lipschitz} holds only on a neighborhood of the reference matrix. To
compare the reference score with that of an arbitrary Hurwitz matrix
differing in the positions $\mathcal{E}$, we combine it with the
boundedness of the simplex outside this neighborhood.

\begin{theorem}[Reference-dependent score robustness]\label{thm:graph}
Fix a Hurwitz matrix $A$ and a criterion $\CS$. There exists $K_A^{\CS}>0$ such that, for every Hurwitz $\widetilde A$ agreeing with $A$ outside $\mathcal E$,
\begin{equation}\label{eq:reference-bound}
 |p_i^{\CS}(A)-p_i^{\CS}(\widetilde A)|
 \leq K_A^{\CS}d_{\mathcal E}(A,\widetilde A)\quad (i=1,\ldots,n).
\end{equation}
\end{theorem}
\begin{IEEEproof}
Choose $\delta >0$ so that $\lVert\widetilde A-A\rVert_{\mathrm F}<\delta$ lies in the neighborhood of Lemma~\ref{lem:lipschitz}. In that ball, \eqref{eq:reference-bound} holds with $L_{\CS}$. Outside the ball, $d_{\mathcal E}(A,\widetilde A)\geq \delta$ and $|p_i^{\CS}(A)-p_i^{\CS}(\widetilde A)|\leq1$. Taking $K_A^{\CS}=\max\{L_{\CS},1/\delta\}$ proves the claim.
\end{IEEEproof}

Theorem~\ref{thm:graph} provides the robustness guarantee motivated in Section~\ref{sec:introduction}. Because $d_{\mathcal{E}}$ allows
zero weights, Theorem~\ref{thm:graph} covers both inaccurate weights and changed connections, including the insertion
and deletion of edges, and is not restricted to small perturbations.
The bound holds for each node separately, so no component~$p_i^{\rm CS}$
can change disproportionately to the model change. In particular, if
two nodes satisfy
$p^{\mathrm{CS}}_i(A) - p^{\mathrm{CS}}_j(A) > 2K^{\mathrm{CS}}_A\, d_{\mathcal{E}}(A,\widetilde A)$,
then node~$i$ remains ranked above node~$j$ at $\widetilde A$. Thus, sufficiently large differences between node scores at the reference matrix guarantee that their ranking is preserved under small model perturbations.

However, the bound in \eqref{eq:reference-bound} can be conservative for a particular model change because the constant $K^{\mathrm{CS}}_A$ must hold uniformly over
all admissible perturbations: it is built from worst-case estimates
over a neighborhood of $A$ and from the trivial bound outside it.
Section~\ref{sec:sensitivity} instead fixes a perturbation direction at the
reference matrix and characterizes the first-order change of each
score exactly, rather than bounding it.

\begin{remark}[Relation to centrality robustness]
Theorem~\ref{thm:graph} gives a reference-dependent robustness bound analogous to centrality robustness proposed in~\cite[Def.~1]{segarra2016stability}. Whereas \cite{segarra2016stability} considers perturbations of positive edge weights on a fixed labeled support, our distance~\eqref{eq:graph-distance} allows signed or zero coefficients and system-matrix diagonal changes within the Hurwitz class.
\end{remark}

\section{Direction-specific local robustness}
\label{sec:sensitivity}
In Section~\ref{sec:continuity}, we bounded score changes uniformly over all
admissible perturbations $D \in \R^{n\times n}$ in~\eqref{eq:matrix-path}. In this section, we fix a perturbation
direction~$D$ at the reference matrix~$A$ and characterize the first-order
change of the scores along this direction, which we call the
\emph{controllability score response vector} for system~\eqref{eq:autonomous}. For single-edge
weakening, its norm provides a robustness-based measure of edge
importance.

\subsection{Characterization and interpretation}

To define the controllability score response vector for system~\eqref{eq:autonomous}, we consider the model-change paths defined in Section~\ref{sec:formulation}, with direction $D$ as in~\eqref{eq:matrix-path}. Fix a criterion \(\mathrm{CS}\in\{\mathrm V,\mathrm A\}\) and a reference matrix $A$. Throughout this section, write
\[
p=p^{\mathrm{CS}}(A),\qquad
g=\nabla_p f_{\mathrm{CS}}(p,A),\qquad
H=\nabla_{pp}^2f_{\mathrm{CS}}(p,A).
\]
The simplex KKT conditions imply that there exist multipliers \(\nu\in\mathbb R\) and \(\mu\in\mathbb R_{\geq0}^n\) such that
\begin{equation}\label{eq:kkt}
 g+\nu\one-\mu=0,\qquad\mu\geq0,\qquad\mu_ip_i=0.
\end{equation}
Define the support and active index sets
\begin{equation}\label{eq:index_set}
 \begin{aligned}
 J&=\{i:p_i>0\},\\
 I_+&=\{i:p_i=0,\mu_i>0\},\\
 I_0&=\{i:p_i=0,\mu_i=0\}.
 \end{aligned}
\end{equation}
The tangent cone and critical cone are
\begin{equation}
\label{eq:critical-cone}
\begin{aligned}
T_\Delta(p)
&= \{z:\mathbf{1}^\top z=0,\;
       z_i\geq 0\text{ if }p_i=0\},\\
\mathcal K
&=T_\Delta(p)\cap g^\perp\\
&=\{z:\mathbf{1}^\top z=0,\;
       z_{I_+}=0,\;z_{I_0}\geq0\}.
\end{aligned}
\end{equation}

The following result characterizes this local response, using the constrained perturbation-analysis principle in~\cite{bonnans2000perturbation}.

\begin{theorem}[Directional derivative of the optimal allocation]\label{thm:derivative}
For the path \eqref{eq:matrix-path}, the right derivative of $p^{\CS}(A(\theta))$ exists at $\theta=0$, and
\begin{equation}\label{eq:score-expansion}
 p^{\CS}(A(\theta))=p+\theta d+o(\theta),
\end{equation}
\begin{equation}\label{eq:response-qp}
 d=\argmin_{z\in\mathcal K}\left\{\tfrac12z^\top Hz+r^\top z\right\},
\end{equation}
\begin{equation}\label{eq:forcing-definition}
 r=r^{\mathrm{CS}}(D):=\left.\partial_\theta\nabla_p f_{\CS}(p,A(\theta))\right|_{\theta=0}.
\end{equation}
\end{theorem}
\begin{IEEEproof}
See Appendix~\ref{app:derivative}.
\end{IEEEproof}

We call the vector $d$ in~\eqref{eq:response-qp} the controllability score response vector for system~\eqref{eq:autonomous}
along the prescribed direction $D$. The
quadratic program~\eqref{eq:response-qp} has a unique solution: its
feasible set $\mathcal K$ is a closed convex cone in $\one^\perp$, on
which $H$ is positive definite by Lemma~\ref{lem:curvature}. Thus $d$
is a unique well-defined vector determined by the reference matrix,
the criterion, and the direction, which makes comparisons of
$\lVert d\rVert_2$ across edges unambiguous.

In~\eqref{eq:score-expansion}, $d_i$ gives the first-order change in the score of node $i$: its sign indicates an increase or decrease. With the perturbation direction and parameter scale fixed, smaller $|d_i|$ means greater first-order local robustness of that node's score along that direction. If $d_i=0$, only the first-order term vanishes. The norm \(\|d\|_2\) measures the first-order rate of change of the entire score vector. Its magnitude depends on the perturbation direction \(D\) and the parameterization of the path.

\begin{example}[Fixed-entry deletion]
Consider the single-edge weakening path~\eqref{eq:deletion}, with \(D=-A_{rs}E_{rs}\), and let \(d\) be the response vector given by Theorem~\ref{thm:derivative}. Then,
\begin{equation}\label{eq:norm-expansion}
 \lVert p^{\CS}(A(\theta))-p\rVert_2
 =\theta \lVert d\rVert_2+o(\theta).
\end{equation}
Appendix~\ref{app:comparison} provides a two-node example showing that neither criterion has uniformly larger response norms.
\end{example}

The weakening of a single edge $s\to r$ turns the response norm $\lVert d\rVert_2$ into a measure of edge importance from the viewpoint of robustness. By~\eqref{eq:norm-expansion}, $\lVert d\rVert_2$ is the first-order rate of change of the optimal allocation under fractional edge weakening. A large $\lVert d\rVert_2$ therefore means that the
controllability score~$p^{\rm CS}(A)$ is sensitive to the edge~$s\to r$: an error in its
estimated weight, or a change in the connection itself, shifts the
allocation more than the same relative change of an edge with a small
$\lVert d\rVert_2$. In this sense, edges with large $\lVert d\rVert_2$ are
important for the robustness of the score.

\begin{remark}[Relation to edge centralities]\label{rem:edge-centrality}
Our robustness-based edge importance differs from existing approaches by measuring how the optimal input allocation responds to edge perturbations. Previous studies assess edge importance through structural
controllability under link removal~\cite{liu2011controllability},
Gramian-based performance under fixed
actuation~\cite{chanekar2021energy,chanekar2022encoding}, or stability and transfer-function responses to edge modifications~\cite{lindmark2023investigating}. The closest to
the present setting is the edge centrality matrix
of~\cite{chanekar2022encoding}, which encodes the first-order change of
a Gramian-based metric when the input matrix is kept fixed. In contrast, our approach examines how the optimal allocation \(p^{\mathrm{CS}}(A)\) changes as \(A\) varies. Because the optimal
allocation is the minimizer of a constrained problem~\eqref{eq:score}, its first-order
change is not obtained by differentiating a metric at a fixed allocation~$p$; it
requires the critical-cone QP of Theorem~\ref{thm:derivative}. Thus
$\lVert d\rVert_2$ measures the importance of an edge for the robustness of
the optimal allocation, rather than for a controllability metric under
a given allocation.
\end{remark}

\begin{remark}[Comparison with reoptimization]\label{rem:reoptimization}
The change of the scores along a path could also be examined by
reoptimizing the allocation at perturbed matrices and comparing the
results with the reference scores. The response vector $d$ avoids two
limitations of this direct approach. First, reoptimization
solves~\eqref{eq:score} again for each perturbed matrix, so comparing
many candidate changes, such as the weakening of every edge, requires
as many reoptimizations; the response vector reuses the reference
allocation and its Gramians and requires one convex QP per direction. Second, a change computed at a finite
perturbation size mixes first-order and higher-order effects and
depends on the chosen size, whereas difference quotients at small sizes
amplify the optimization tolerance; the response vector depends only
on the reference matrix, the criterion, and the direction.
\end{remark}

\subsection{Computation of the response vector}

To compute the response vector $d$, we need the derivative $r$ of the gradient
$\nabla_p f_{\CS}$ with respect to the prescribed change of the model
in~\eqref{eq:forcing-definition}. The next lemma evaluates it through
differentiated Lyapunov equations, with $p=p^{\CS}(A)$ held fixed.

\begin{lemma}[Gradient derivative via Lyapunov equations]\label{lem:forcing}
At fixed $p$, let $R=W_A(p)^{-1}$. For the direction $D$ in~\eqref{eq:matrix-path}, define $Y_i[D]$ by
\begin{equation}\label{eq:forcing-lyapunov}
 AY_i[D]+Y_i[D]A^\top+DW_i+W_iD^\top=0,
\end{equation}
and set $Y[D]=\sum_i p_iY_i[D]$. The gradient derivative $r=r^{\CS}(D)$ in~\eqref{eq:forcing-definition} has components
\begin{equation}\label{eq:forcing-V}
 r_i^{\mathrm V}(D)=\tr(RY[D]RW_i)-\tr(RY_i[D]),
\end{equation}
\begin{equation}\label{eq:forcing-A}
 r_i^{\mathrm A}(D)=2\tr(RY[D]RW_iR)-\tr(RY_i[D]R).
\end{equation}
The gradient derivative $r^{\mathrm{CS}}(D)$ is linear in $D$. In particular,
\[
 r^{\CS}(D)=\sum_{r,s}D_{rs}r^{\mathrm{CS}}(E_{rs}).
\]
\end{lemma}
\begin{IEEEproof}
All derivatives below are taken with respect to \(\theta\) at \(\theta=0\) along the path \(A(\theta)\), with \(p\) held fixed. Differentiating \eqref{eq:lyapunov} shows that the derivative of \(W_i(A(\theta))\) satisfies \eqref{eq:forcing-lyapunov}. By uniqueness of the Lyapunov solution, this derivative equals \(Y_i[D]\). The derivative of $R$ in direction $D$ is $R'[D] = -R\,Y[D]\,R$. Differentiating \eqref{eq:gradients} yields \eqref{eq:forcing-V} and \eqref{eq:forcing-A}; in the latter, the two positive trace terms coincide because $R$, $Y[D]$, and $W_i$ are symmetric and $\tr(M)=\tr(M^\top)$. The Lyapunov solution and the partial gradient derivative are linear in $D$.
\end{IEEEproof}

Therefore, we can compute the response vector by using Algorithm~\ref{alg:response}, which combines the reference optimization, differentiated Lyapunov equations, and critical-cone QP. Once the reference allocation is available, this procedure gives the first-order response for each specified direction without reoptimizing at finely spaced perturbation levels.
\begin{algorithm}[Controllability score response vector]\label{alg:response}
\begin{algorithmic}[1]
\Require Hurwitz reference matrix $A$, choice $\CS\in\{\mathrm V,\mathrm A\}$, direction $D$.
\Ensure Controllability score response vector $d$.
\State Solve~\eqref{eq:lyapunov} for $W_i$, $i=1,\ldots,n$.\label{line:gramians}
\State Solve~\eqref{eq:score} for $p=p^{\CS}(A)$.\label{line:baseline}
\State Compute \(g=\nabla_p f_{\mathrm{CS}}(p,A)\) and \(H=\nabla_{pp}^2f_{\mathrm{CS}}(p,A)\) using~\eqref{eq:gradients}, \eqref{eq:HV}, \eqref{eq:HA}.\label{line:gradient}\label{line:hessian}
\State Set $\nu=-g_j$ for any $j$ with $p_j>0$ and $\mu=g+\nu\one$.\label{line:nu}\label{line:mu}
\State Form $I_+=\{i:p_i=0,\ \mu_i>0\}$ and $I_0=\{i:p_i=0,\ \mu_i=0\}$.\label{line:ipositive}\label{line:izero}
\State Form $\mathcal K$ using~\eqref{eq:critical-cone}.\label{line:cone}
\State Solve~\eqref{eq:forcing-lyapunov} for $Y_i[D]$.\label{line:direct-lyapunov}
\State Form $Y[D]=\sum_i p_iY_i[D]$.\label{line:aggregate}
\State Evaluate $r$ using~\eqref{eq:forcing-V} or~\eqref{eq:forcing-A}.\label{line:forcing}
\State Solve the QP~\eqref{eq:response-qp} for $d$.\label{line:qp}
\State \Return $d$.
\end{algorithmic}
\end{algorithm}

The computation uses standard Lyapunov and convex optimization solvers. In lines~\ref{line:gramians} and~\ref{line:direct-lyapunov} of Algorithm~\ref{alg:response}, equations~\eqref{eq:lyapunov} and~\eqref{eq:forcing-lyapunov} are continuous-time Lyapunov equations and can be solved by the Bartels--Stewart method~\cite{bartels1972solution}. In line~\ref{line:baseline}, the reference allocation~\eqref{eq:score} can be computed by the projected-gradient method with an Armijo line search described in~\cite[Sec.~V, Algorithms~1--2 and Remark~5]{sato2024scores}, using the VCS or AECS objective and its gradient.

The response problem in line~\ref{line:qp} is a convex quadratic program over the critical cone \(K\), which can be expressed using linear equality and inequality constraints. Specifically, the fixed components \(z_{I_+}=0\) can be eliminated, while the equality \(\one^\top z=0\) and nonnegativity constraints \(z_{I_0}\geq0\) are retained. The resulting QP can be solved, for example, by a primal-dual interior-point method~\cite[Secs.~4.4 and~11.7]{boyd2004convex}.

Under strict complementarity (\(I_0=\emptyset\)), the critical-cone QP reduces to a linear system, as shown in Appendix~\ref{app:derivative}. In this case, the critical cone contains only equality constraints, since the inactive allocation components are fixed at zero. The resulting bordered linear system has a coefficient matrix independent of \(D\), so its factorization can be reused for different perturbation directions, replacing repeated QP solves with linear solves.

\begin{remark}[Adjoint evaluation for multiple directions]
\label{rem:adjoint-forcing}
For multiple directions at the same reference, replacing lines~\ref{line:direct-lyapunov}--\ref{line:forcing} of Algorithm~\ref{alg:response} by an adjoint evaluation makes evaluation of $r$ more efficient by reusing direction-independent solutions. At fixed $A$ and $p$, set $W=W_A(p)$ and $R=W^{-1}$. For VCS, let $M_0=R$ and $M_i=RW_iR$; for AECS, let $M_0=R^2$ and $M_i=RW_iR^2+R^2W_iR$. Solve
\[
 A^\top P_j+P_jA+M_j=0,\qquad j=0,\ldots,n.
\]
Lyapunov adjointness gives $\langle M_j,Y[D]\rangle_{\rm F}=2\langle P_jW,D\rangle_{\rm F}$, so the gradient derivative in Lemma~\ref{lem:forcing} can be evaluated as
\[
 r_i^{\CS}(D)=2\langle P_iW-P_0W_i,D\rangle_{\rm F},
\]
where $\langle U,V\rangle_{\rm F}=\tr(U^\top V)$. The coefficient matrices \(2(P_iW-P_0W_i)\) are independent of \(D\) and can be reused to evaluate the gradient derivative for multiple perturbation directions.

The adjoint method reduces the cost of evaluating the gradient derivative for multiple perturbation directions. After an \(O(n^4)\) preprocessing step, each evaluation of \(r_{\mathrm{CS}}(D)\) requires only \(O(nk)\) operations when \(D\) has \(k\) nonzero entries, compared with \(O(n^4)\) operations per direction using the direct method. In particular, for single-edge weakening, \(k=1\) or \(2\), so the cost is reduced to \(O(n)\) per direction. The QP for \(d\) must still be solved separately for each perturbation direction.
\end{remark}

\begin{figure*}[t]
\centering
\includegraphics[width=0.85\textwidth]{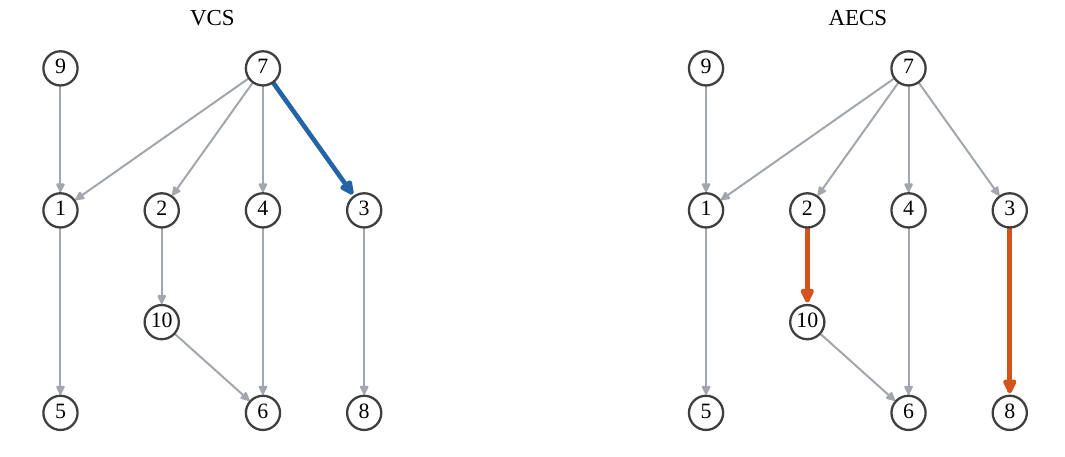}
\caption{Ten-node directed network used for the edge-weakening experiments. Each arrow has weight $0.2$; the two state matrices are defined in Section~\ref{sec:numerical-setup}. The blue edge and the two orange edges mark the largest VCS response norm and the two largest AECS response norms, respectively, under Laplacian-preserving weakening of $A_{\rm L}$. See Table~\ref{tab:response-norms} for the values.}
\label{fig:edge-response}
\end{figure*}

\section{Numerical experiments}\label{sec:numerical}\label{sec:numerical-setup}

We examine how edge weakening redistributes optimal allocations and how accurately the response vector predicts finite changes. A ten-node example illustrates the response patterns, and 100-node networks test prediction accuracy over a range of perturbation sizes.

The experiments in Sections~\ref{sec:numerical-setupA}
and~\ref{sec:numerical-setupB} use the following common setup. For each
directed graph with $n$ nodes, let $\operatorname{Adj}_{uv}=0.2$ if
$u\to v$ is an edge and zero otherwise, and let $L$ be the row
Laplacian of Example~\ref{ex:laplacian} built from this
$\operatorname{Adj}$. We consider the two dynamical models
\[
 A_{\rm L}=-L-0.002I,\qquad A_{\rm F}=\operatorname{Adj}^{\top}-I.
\]
The model $A_{\rm L}$ is the row-Laplacian model of
Example~\ref{ex:laplacian} with $\alpha=0.002$, and both models are
Hurwitz. For relative weakening of an edge $s\to r$, we use
\[
 A_\bullet(\theta)=A_\bullet+\theta D_\bullet,\qquad
 \bullet\in\{\mathrm L,\mathrm F\},
\]
where
\[
 D_{\rm L}=0.2(E_{rr}-E_{rs}),\qquad D_{\rm F}=-0.2E_{rs}.
\]
Thus $D_{\rm L}$ is the Laplacian-preserving path of
Example~\ref{ex:laplacian} with $w=0.2$, and $D_{\rm F}$ is the
fixed-entry path of Example~\ref{ex:edge_weakening} with $a=0.2$. Here,
$\theta=0$ corresponds to the original edge weight, whereas $\theta=1$
removes the edge. Both paths remain Hurwitz for $0\leq\theta\leq1$:
for $A_{\rm L}$, $L(\theta)$ remains a row Laplacian with nonnegative
edge weights, and $A_{\rm F}$ is Metzler, so
Theorem~\ref{thm:metzler} in Appendix~\ref{app:metzler} applies. The controllability score response
vector $d$ is computed by Algorithm~\ref{alg:response}.

\subsection{Allocation responses to edge weakening}\label{sec:numerical-setupA}

\begin{figure*}[t]
\centering
\includegraphics[width=0.98\textwidth]{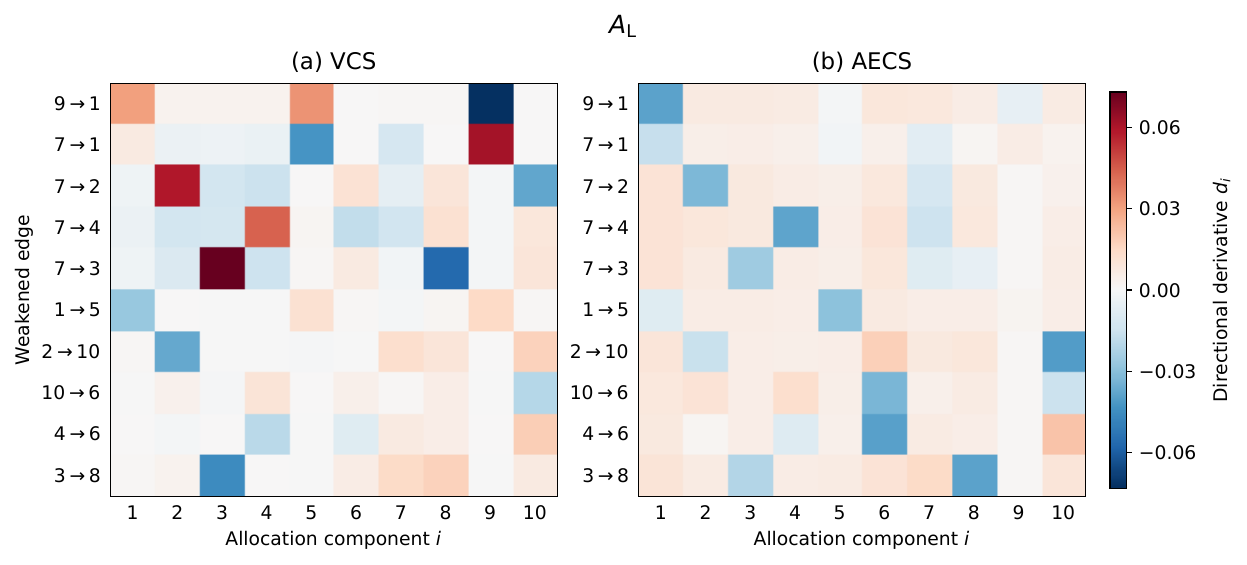}
\caption{Node-wise response to Laplacian-preserving edge weakening for $A_{\rm L}$. Each row gives the response vector $d$ for the indicated edge; each column identifies the node whose allocation changes. The VCS (left) and AECS (right) panels share a linear color scale centered at zero: red indicates an increase and blue a decrease per unit fractional weakening. The range is approximately $\pm0.0731$; Fig.~\ref{fig:components-F} uses a different range.}
\label{fig:components-L}
\end{figure*}

\begin{figure*}[t]
\centering
\includegraphics[width=0.98\textwidth]{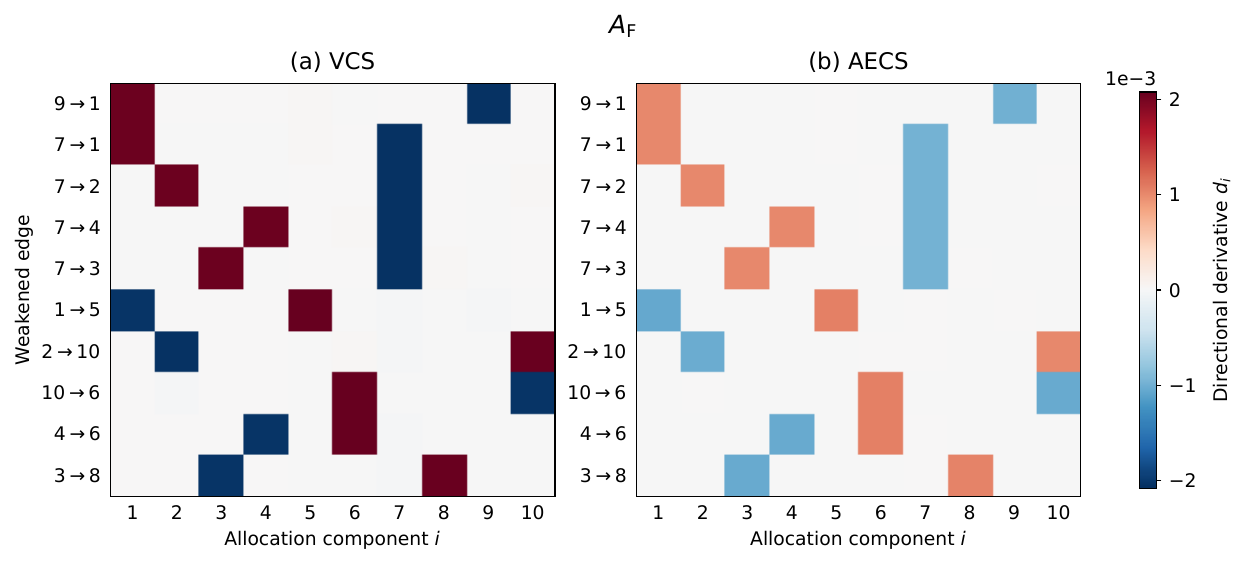}
\caption{Node-wise response to fixed-diagonal edge weakening for $A_{\rm F}$. Each row gives the response vector $d$ for the indicated edge; each column identifies the node whose allocation changes. The VCS (left) and AECS (right) panels share a linear color scale centered at zero: red indicates an increase and blue a decrease per unit fractional weakening. The range is approximately $\pm0.00208$; Fig.~\ref{fig:components-L} uses a different range.}
\label{fig:components-F}
\end{figure*}

\begin{table}[t]
\centering
\caption{Response norms $\lVert d\rVert_2$ per unit fractional weakening of each edge in the ten-node network. L denotes $A_{\rm L}$ with Laplacian-preserving weakening, and F denotes $A_{\rm F}$ with fixed-diagonal weakening. Rank emphasis is limited to the two L columns; the F columns are not ranked because several displayed values are tied or nearly tied.}
\label{tab:response-norms}
\small
\begin{tabular}{@{}lrrrr@{}}
\hline
Edge & VCS L & AECS L & VCS F & AECS F\\
\hline
$9\to1$ & \textcolor{RankSecond}{\underline{0.0857}} & 0.0436 & 0.002899 & 0.001420\\
$7\to1$ & \textcolor{RankThird}{\textit{0.0761}} & 0.0214 & 0.002892 & 0.001405\\
$7\to2$ & 0.0746 & 0.0390 & 0.002898 & 0.001402\\
$7\to4$ & 0.0541 & 0.0463 & 0.002897 & 0.001401\\
$7\to3$ & \textcolor{RankFirst}{\textbf{0.0949}} & 0.0325 & 0.002893 & 0.001402\\
$1\to5$ & 0.0331 & 0.0346 & 0.002915 & 0.001494\\
$2\to10$ & 0.0438 & \textcolor{RankSecond}{\underline{0.0500}} & 0.002915 & 0.001440\\
$10\to6$ & 0.0240 & 0.0422 & 0.002905 & 0.001476\\
$4\to6$ & 0.0293 & \textcolor{RankThird}{\textit{0.0475}} & 0.002905 & 0.001474\\
$3\to8$ & 0.0515 & \textcolor{RankFirst}{\textbf{0.0506}} & 0.002901 & 0.001468\\
\hline
\end{tabular}
\par\smallskip
\scriptsize Top three in each L column: \textcolor{RankFirst}{\textbf{1st (bold)}}, \textcolor{RankSecond}{\underline{2nd (underlined)}}, \textcolor{RankThird}{\textit{3rd (italic)}}.
\end{table}

For the ten-node graph in Fig. 2, Table I shows that edge importance depends on the scoring criterion and the dynamical model. Under \(A_{\rm L}\), the largest VCS response norms occur on edges leaving source nodes, whereas the largest AECS response norms occur on edges farther downstream.

Figs.~\ref{fig:components-L}--\ref{fig:components-F} show that the node-wise response patterns also differ between the models. Under \(A_{\rm F}\), the response is concentrated mainly on the sender \(s\) and receiver \(r\), with \(d_s<0<d_r\). Under \(A_{\rm L}\), the response extends to other nodes, and some components have opposite signs for VCS and AECS.

These differences indicate that the response is not determined by the
graph topology alone. For the same graph and the same weakened edge,
$A_{\rm L}$ and $A_{\rm F}$ give different response norms, and some
components of the response vectors have opposite signs. Previous
studies have assessed edge importance from the perspective of
structural controllability, including the effects of link failures and
edge removal~\cite{rahimian2013structural,mengiste2015edge,chung2021dilations}. These approaches focus on the zero--nonzero patterns of the system and input matrices. In contrast, our approach quantifies changes in the optimal allocation under specified edge-weight perturbations. Our results show that this response depends on the specified dynamical model and perturbation path, even when the underlying graph is unchanged. Therefore it is useful
to report them together with the edge importance.

\subsection{Accuracy of the first-order prediction}\label{sec:numerical-setupB}

We test nine directed graphs with 100 nodes and 200 edges:
three realizations each of the Erd\H{o}s--R\'enyi
$G(n,m)$ (GNM) random graph model~\cite{erdos1959random},
the Barab\'asi--Albert (BA) preferential-attachment
model~\cite{barabasi1999emergence}, and the
Watts--Strogatz (WS) small-world
model~\cite{watts1998collective}. For each graph, we randomly select three distinct edges before optimization. We compare the response vector $d$ with reoptimized allocations $p(\theta)=p^{\CS}(A_\bullet(\theta))$ at $\theta\in\{0.001,0.01,0.03,0.1,0.3,1\}$ through
\[
 \varepsilon(\theta)=\frac{\|(p(\theta)-p(0))/\theta-d\|_2}
 {\max\{\|d\|_2,\epsilon_{\rm mach}\}},
\]
where $\epsilon_{\rm mach}$ is double-precision machine epsilon.

\begin{table*}[t]
\centering\small
\caption{Relative derivative error $\varepsilon(\theta)$ in percent on the nine 100-node networks. Each cell reports the median (maximum) over 27 graph--edge cases for the indicated model and criterion, rounded to three decimal places. The same nine graphs and three selected edges per graph are used at every fractional weakening $\theta$; $\theta=1$ corresponds to full edge deletion.}
\label{tab:hundred-width}
\setlength{\tabcolsep}{10pt}
\begin{tabular}{@{}crrrr@{}}
\hline
$\theta$ & $A_{\rm L}$/VCS & $A_{\rm L}$/AECS & $A_{\rm F}$/VCS & $A_{\rm F}$/AECS\\
\hline
$0.001$ & 0.041 (0.103) & 0.014 (0.078) & 0.052 (0.052) & 0.052 (0.057)\\
$0.01$ & 0.409 (1.031) & 0.141 (0.779) & 0.516 (0.523) & 0.519 (0.574)\\
$0.03$ & 1.241 (3.130) & 0.427 (2.335) & 1.546 (1.569) & 1.555 (1.721)\\
$0.1$ & 4.299 (10.901) & 1.455 (7.673) & 5.139 (5.218) & 5.171 (5.727)\\
$0.3$ & 14.165 (37.409) & 4.679 (20.951) & 15.313 (15.547) & 15.413 (17.114)\\
$1$ & 72.795 (650.101) & 21.828 (111.581) & 50.246 (50.916) & 50.484 (56.519)\\
\hline
\end{tabular}
\end{table*}

Table~\ref{tab:hundred-width} shows median errors below 0.6\% for every model–criterion combination at $\theta\leq0.01$. Errors grow at larger $\theta$; at full edge deletion, medians range from 21.8\% to 72.8\%, with a 650.1\% maximum for $A_{\rm L}$/VCS. These fixed-sample results support the first-order prediction for small perturbations while illustrating the limitations of a local approximation.

\section{Conclusion}
\label{sec:conclusion}
This paper analyzed the robustness of VCS and AECS to
perturbations of the system matrix. For Hurwitz systems
with node-wise inputs and a simplex allocation constraint,
we established local Lipschitz continuity of both scores,
including at boundary allocations. We then obtained a
reference-dependent robustness bound for every admissible
Hurwitz comparison matrix, beyond small perturbations.
For a prescribed differentiable matrix path, we characterized
the directional derivative of the optimal allocation as the
unique solution of a critical-cone quadratic program,
with the gradient derivative computed from differentiated
Lyapunov equations.

The score response provides a robustness-based measure
of edge importance and reveals its dependence on the
scoring criterion and the specified dynamical model.
For fractional weakening of a single edge, the response
norm $\|d\|_2$ measures the first-order change in the
optimal allocation. The two-node analytical example
in Appendix~\ref{app:comparison} shows that neither
VCS nor AECS has uniformly larger response norms.
Numerical experiments further demonstrate that the
response magnitude and node-wise redistribution depend
on the criterion, dynamical model, and associated
weakening path, even for the same underlying graph topology.

Earlier work introduced the VCS and AECS formulations and analyzed their uniqueness~\cite{sato2024scores,sato2025uniqueness}. Infinite-horizon scoring has also been studied through limiting Gramian constructions~\cite{umezu2026infinite}, and the restricted injectivity of the Gramian map provides quantitative curvature bounds for a generalized allocation formulation~\cite{sato2026projected}. These developments motivate extending the present robustness analysis to finite-horizon scores and to limiting formulations beyond Hurwitz dynamics. Such extensions require checking the relevant uniqueness, curvature, and parameter-dependence properties; they are not consequences of the results proved here.

Future work includes characterizing how network structure and system dynamics determine the score response, as well as separating the effects of the reference model and perturbation direction through controlled comparisons.

\section*{Acknowledgment}
This work was supported by JST PRESTO, Japan, Grant Number JPMJPR25K4.

\appendices

\section{Proof of Lemma~\ref{lem:curvature}}\label{app:curvature}
We prove Lemma~\ref{lem:curvature} by deriving explicit lower bounds
on the Hessian quadratic forms restricted to
$\mathbf{1}^{\perp}$. The key step is to establish
a positive lower bound for the Gramian map
$Q_A$ on this subspace.

For every Hurwitz $A$, the map $Q_A:\mathbb R^n\to\mathbb S^n$ defined in~\eqref{eq:QA} is injective. Let $Z\in\mathbb R^{n\times(n-1)}$ have orthonormal columns that span $\one^\perp$, and define
\begin{align}
 M_A &=[\vecop W_1(A)\ \cdots\ \vecop W_n(A)], \,\,
\kappa_A =\sigma_{\min}(M_AZ)>0,\\
\beta_A &=\max_i\lambda_{\max}(W_i(A)),\,\,
 \rho_A=\kappa_A/\beta_A.
\end{align}

We claim that the following bounds hold for every \(p\in X_A\) and \(z\in\mathbf1^\perp\):
\begin{equation}\label{eq:curvature-V-explicit}
 z^\top H_{\mathrm V}(p,A)z\geq\rho_A^2\lVert z\rVert_2^2,
\end{equation}
\begin{equation}\label{eq:curvature-A-explicit}
 z^\top H_{\mathrm A}(p,A)z\geq\frac{2\rho_A^2}{\beta_A}\lVert z\rVert_2^2.
\end{equation}
\begin{IEEEproof}
Multiply the $i$th equation in \eqref{eq:lyapunov} by $z_i$ and sum over $i$. Linearity of the Lyapunov operator and the identity $\sum_i z_ie_ie_i^\top=\diag(z)$ give
\begin{equation}\label{eq:summed-lyapunov}
 AQ_A(z)+Q_A(z)A^\top+\diag(z)=0.
\end{equation}
If $Q_A(z)=0$, \eqref{eq:summed-lyapunov} reduces to $\diag(z)=0$, hence $z=0$. Thus $Q_A$ is injective on $\mathbb R^n$, and in particular on $\one^\perp$. For $z\in\one^\perp$, write $z=Zy$, where \(y=Z^\top z\). Since the columns of $Z$ are orthonormal, $\|y\|_2=\|z\|_2$, and
\[
 \|Q_A(z)\|_{\mathrm F}=\|M_AZy\|_2
 \geq\kappa_A\|z\|_2.
\]
Since $Q_A$ is injective, $M_AZ$ has full column rank, so $\kappa_A>0$. The singular values of $M_AZ$ do not depend on the choice of $Z$.

Each $W_i(A)$ is positive semidefinite and bounded above by $\beta_AI$. Since $p\in\Delta$, the convex combination~\eqref{eq:WA} satisfies $0\prec W_A(p)\preceq\beta_AI$. Following the argument in~\cite[Lem.~6]{sato2026projected}, we obtain~\eqref{eq:curvature-V-explicit} and~\eqref{eq:curvature-A-explicit}.

Thus, Lemma~\ref{lem:curvature} holds with $m_{\mathrm V}(A)=\rho_A^2$ and $m_{\mathrm A}(A)=2\rho_A^2/\beta_A$.
\end{IEEEproof}

\section{Proof of Lemma~\ref{lem:lipschitz}}\label{app:lipschitz}

We prove local Lipschitz continuity by combining
the curvature bound of Lemma~\ref{lem:curvature} with a perturbation
bound on the gradient. The proof has two steps.
First, we establish a common neighborhood in
which optimal allocations are feasible for both
optimization problems. We then compare their
optimality conditions to bound the difference
between the allocations.

\subsection{Uniform bounds and cross-feasibility}
We first establish a uniform positive-definiteness bound for the optimal Gramians near the reference matrix. This bound is obtained by comparing the optimal objective values with those at the uniform allocation \(p_0\), without assuming continuity of the optimizers.

Take a closed ball $\mathcal{B}$ about $A_0$ contained in the Hurwitz
domain. Since each \(W_i(A)\) depends smoothly on \(A\), both \(\beta_A\) and \(\kappa_A\) are continuous on \(\mathcal B\). Compactness and convexity of
$\mathcal{B}$ therefore give $\bar\beta > 0$, $\kappa > 0$, and
$\ell > 0$ such that for all $A, \widetilde A \in \mathcal{B}$,
\begin{align}
\beta_A &\le \bar\beta, \qquad \kappa_A \ge \kappa, \label{eq:uniform-curvature-data}\\
\max_i \|W_i(A) - W_i(\widetilde A)\|_{\mathrm F} &\le \ell \|A - \widetilde A\|_{\mathrm F}. \label{eq:Wi-lipschitz}
\end{align}

For $p^0=\one/n$, $W_A(p^0)$ is positive definite for every $A$ in the ball. Continuity of $A\mapsto W_A(p^0)$, compactness of $\mathcal B$, and continuity of the determinant and inverse on positive-definite matrices give
\begin{equation}\label{eq:D0}
 D_0=\min_{A\in\mathcal B}\det W_A(p^0)>0,
\end{equation}
\begin{equation}\label{eq:C0}
 C_0=\max_{A\in\mathcal B}\tr\bigl(W_A(p^0)^{-1}\bigr)<\infty.
\end{equation}
At the VCS optimum, a comparison with $p^0$ gives $-\log\det W_A(p^{\mathrm V}(A))\leq-\log\det W_A(p^0)$, so its determinant is at least $D_0$. At the AECS optimum, the corresponding comparison bounds the inverse trace by $C_0$. Consequently,
\begin{equation}\label{eq:eta-V}
 \lambda_{\min}(W_A(p^{\mathrm V}(A)))
 \geq D_0/\bar\beta^{n-1}=:\eta_{\mathrm V},
\end{equation}
\begin{equation}\label{eq:eta-A}
 \lambda_{\min}(W_A(p^{\mathrm A}(A)))\geq C_0^{-1}=:\eta_{\mathrm A}.
\end{equation}
For \eqref{eq:eta-V}, the determinant of the optimal Gramian is at least $D_0$, while each of its other eigenvalues is at most $\bar\beta$. For \eqref{eq:eta-A}, the reciprocal of the least eigenvalue is at most the trace of the inverse, which is at most $C_0$.

Fix the criterion and write $\eta=\eta_{\CS}$. Reduce to an open ball $U$ about $A_0$ whose diameter is at most $\eta/(2\ell)$. Choose $U$ inside $\mathcal B$. For every $p\in\Delta$, the triangle inequality, $p_i\geq0$, $\sum_i p_i=1$, and \eqref{eq:Wi-lipschitz} give
\begin{equation}
 \begin{aligned} 
 \|W_A(p)-W_{\widetilde A}(p)\|_2
 &\leq\sum_i p_i\|W_i(A)-W_i(\widetilde A)\|_2\\
 &\leq\ell\|A-\widetilde A\|_{\mathrm F}. 
 \end{aligned}\label{eq:Gramian-difference}
\end{equation}
For $p_{\widetilde A}=p^{\CS}(\widetilde A)$, \eqref{eq:eta-V} or \eqref{eq:eta-A} gives $W_{\widetilde A}(p_{\widetilde A})\succeq\eta I$. By the eigenvalue perturbation inequality, \eqref{eq:Gramian-difference}, and the diameter bound on $U$, $\lambda_{\min}(W_A(p_{\widetilde A}))\geq\eta-\ell\|A-\widetilde A\|_{\mathrm F}\geq\eta/2$. Thus $p_{\widetilde A}\in X_A$ and
\begin{equation}\label{eq:cross-safe}
 W_A(p^{\CS}(\widetilde A))\succeq(\eta/2)I.
\end{equation}
By interchanging \(A\) and \(\widetilde A\), we also have \(p^{\mathrm{CS}}(A)\in X_{\widetilde A}\). Thus, the optimal allocations at two nearby matrices are feasible for both optimization problems, allowing their optimality conditions to be compared in the next subsection.

\subsection{Gradient perturbation and strong monotonicity}
We now bound the change in the gradient between
two nearby system matrices and combine this estimate
with Lemma~\ref{lem:curvature} and the cross-feasibility established
above to prove Lemma~\ref{lem:lipschitz}.

Set $m=\eta/2$, $p_{\widetilde A}=p^{\CS}(\widetilde A)$, and write $R_A=W_A(p_{\widetilde A})^{-1}$ and $R_{\widetilde A}=W_{\widetilde A}(p_{\widetilde A})^{-1}$. Equation~\eqref{eq:cross-safe} bounds both inverse norms by $m^{-1}$. The inverse identity
\[
 R_A-R_{\widetilde A}=R_A\bigl(W_{\widetilde A}(p_{\widetilde A})-W_A(p_{\widetilde A})\bigr)R_{\widetilde A}
\]
and \eqref{eq:Gramian-difference} give the first estimate in the following. For the second, use the noncommutative identity $R_A^2-R_{\widetilde A}^2=R_A(R_A-R_{\widetilde A})+(R_A-R_{\widetilde A})R_{\widetilde A}$ and submultiplicativity:
\[
 \begin{aligned}
 \lVert R_A-R_{\widetilde A}\rVert_2&\leq\ell m^{-2}\lVert A-\widetilde A\rVert_{\mathrm F},\\
 \lVert R_A^2-R_{\widetilde A}^2\rVert_2&\leq2\ell m^{-3}\lVert A-\widetilde A\rVert_{\mathrm F}.
 \end{aligned}
\]

For VCS, we have
\[
 \begin{aligned}
 &\tr(R_AW_i(A)) -\tr(R_{\widetilde A}W_i(\widetilde A))\\
 &=\tr((R_A-R_{\widetilde A})W_i(A))
   +\tr(R_{\widetilde A}(W_i(A)-W_i(\widetilde A))).
 \end{aligned}
\]
Since $W_i(A)$ is positive semidefinite, the absolute value of the first term is bounded by $\|R_A-R_{\widetilde A}\|_2\tr W_i(A)$, with $\tr W_i(A)\leq n\bar\beta$. The second term is bounded by the Frobenius Cauchy--Schwarz inequality and $\|R_{\widetilde A}\|_{\mathrm F}\leq\sqrt n/m$. For AECS, the cyclicity of the trace rewrites its gradient as $-\tr(R^2W_i)$; the identical splitting uses the second inverse estimate and $\|R_{\widetilde A}^2\|_{\mathrm F}\leq\sqrt n/m^2$. Thus, each gradient component is bounded by the corresponding parentheses in \eqref{eq:CV} or \eqref{eq:CA}, multiplied by $\ell\|A-\widetilde A\|_{\mathrm F}$. Passing from componentwise bounds to the Euclidean norm introduces the factor $\sqrt n$ and gives
\begin{equation}\label{eq:gradient-perturbation}
 \begin{split}
 &\lVert\nabla_p f_{\CS}(p_{\widetilde A},A)-\nabla_p f_{\CS}(p_{\widetilde A},\widetilde A)\rVert_2\\
 &\hspace{4em}\leq C_{\CS}\lVert A-\widetilde A\rVert_{\mathrm F},
 \end{split}
\end{equation}
where
\begin{equation}\label{eq:CV}
 C_{\mathrm V}=\sqrt n\,\ell\bigl(n\bar\beta m^{-2}+\sqrt n\,m^{-1}\bigr),
\end{equation}
\begin{equation}\label{eq:CA}
 C_{\mathrm A}=\sqrt n\,\ell\bigl(2n\bar\beta m^{-3}+\sqrt n\,m^{-2}\bigr).
\end{equation}

Let $p_A=p^{\CS}(A)$ and $\delta p=p_{\widetilde A}-p_A$. Along the segment between the two allocations at fixed $A$, feasibility follows from convexity and \eqref{eq:cross-safe}. Moreover, $\one^\top\delta p=0$. Integrating the Hessian along this segment gives
\[
 \begin{aligned}
 &[\nabla_p f_{\CS}(p_{\widetilde A},A)-\nabla_p f_{\CS}(p_A,A)]^\top\delta p\\
 &=\int_0^1\delta p^\top
 H_{\CS}(p_A+t\delta p,A)\delta p\,dt.
 \end{aligned}
\]
Therefore, Lemma~\ref{lem:curvature} and \eqref{eq:uniform-curvature-data} imply
\begin{equation}\label{eq:strong-monotonicity}
 [\nabla_p f_{\CS}(p_{\widetilde A},A)-\nabla_p f_{\CS}(p_A,A)]^\top\delta p\geq\underline m_{\CS}\lVert\delta p\rVert_2^2,
 \end{equation}
with $\underline m_{\mathrm V}=\kappa^2/\bar\beta^2$ and $\underline m_{\mathrm A}=2\kappa^2/\bar\beta^3$. Optimality at $A$ can be tested in the direction $p_{\widetilde A}-p_A$, because the segment is feasible. Similarly, cross-feasibility with $A$ and $\widetilde A$ interchanged allows optimality at $\widetilde A$ to be tested toward $p_A$. Taking the one-sided directional derivative of each objective along these feasible segments yields
\[
 \nabla_p f_{\CS}(p_A,A)^\top\delta p\geq0,\qquad
 \nabla_p f_{\CS}(p_{\widetilde A},\widetilde A)^\top\delta p\leq0.
\]
Subtracting these inequalities and using \eqref{eq:strong-monotonicity}, we obtain
\[
 \begin{aligned}
 \underline m_{\CS}\|\delta p\|_2^2
 &\leq[\nabla_p f_{\CS}(p_{\widetilde A},A)
       -\nabla_p f_{\CS}(p_{\widetilde A},\widetilde A)]^\top\delta p\\
 &\leq C_{\CS}\|A-\widetilde A\|_{\mathrm F}\|\delta p\|_2,
 \end{aligned}
\]
where the last step uses Cauchy--Schwarz and \eqref{eq:gradient-perturbation}. If $\delta p=0$, the desired estimate is immediate. Otherwise, divide by $\underline m_{\CS}\|\delta p\|_2$ to obtain \eqref{eq:lipschitz} with $L_{\CS}=C_{\CS}/\underline m_{\CS}$. \hfill$\square$

\section{Proof of Theorem~\ref{thm:derivative} and Fixed-Face Reduction}
\label{app:derivative}
We prove Theorem~\ref{thm:derivative} by showing that the difference quotients of the optimal allocation converge to the unique solution of the critical-cone quadratic program~\eqref{eq:response-qp}.

Let $A(\theta)$ satisfy~\eqref{eq:matrix-path}, and write $h_\theta(u)=f_{\CS}(u,A(\theta))$ and $p_\theta=p^{\CS}(A(\theta))$. Set
\[
d_\theta=\frac{p_\theta-p}{\theta},
\qquad \theta>0.
\]
We first show that $d_\theta$ is bounded and
that every subsequential limit belongs to
the critical cone $K$ defined in~\eqref{eq:critical-cone}.

The path satisfies $A(\theta)-A=\theta D+o(\theta)$, hence $\|A(\theta)-A\|_{\mathrm F}=O(\theta)$. For small $\theta\geq0$, it remains in the Hurwitz neighborhood of Lemma~\ref{lem:lipschitz}. Applying Lemma~\ref{lem:lipschitz} with $A_0=A$ and
$\widetilde A=A(\theta)$ gives $\|p_\theta-p\|_2=O(\theta)$. Thus, $d_\theta$ is bounded for sufficiently small $\theta>0$. Since $p_\theta,p\in\Delta$, we have $\one^\top d_\theta=0$. If $p_i=0$, then $(d_\theta)_i=(p_\theta)_i/\theta\geq0$. Therefore, by the definition of $T_\Delta(p)$ in~\eqref{eq:critical-cone}, every $d_\theta$ belongs to $T_\Delta(p)$, and so does every cluster point because $T_\Delta(p)$ is closed.

Using the base KKT conditions \eqref{eq:kkt}, we obtain
\[
 g^\top d_\theta=(\mu-\nu\one)^\top d_\theta
 =\sum_{i:p_i=0}\mu_i(d_\theta)_i\geq0.
\]
Here complementarity eliminates the terms with $p_i>0$, and the remaining factors are nonnegative. By continuity, $W_{A(\theta)}(p)$ is positive definite for small $\theta$. Thus, both $p$ and $p_\theta$ belong to
$X_{A(\theta)}$. Since this set is convex,
the segment joining them also lies in
$X_{A(\theta)}$. By the first-order optimality condition at $p_\theta$,
\[
\nabla_u h_\theta(p_\theta)^\top(p-p_\theta)\ge0.
\]
Since $p-p_\theta=-\theta d_\theta$, this gives
\[
\nabla_u h_\theta(p_\theta)^\top d_\theta\le0.
\]
A first-order expansion of $\nabla_p f_{\mathrm{CS}}$ in $(u,A)$ in a
common positive-definite neighborhood, together with \eqref{eq:forcing-definition}, gives
\begin{equation}
   \nabla_u h_\theta(p_\theta) = g + \theta(Hd_\theta + r) + o(\theta).
\end{equation}
Taking the inner product with the bounded $d_\theta$ in the perturbed inequality gives
\[
 0\leq g^\top d_\theta
 \leq-\theta d_\theta^\top(Hd_\theta+r)+o(\theta)
 =O(\theta).
\]
Along any convergent subsequence $d_{\theta_k}\to\bar d$, it follows that $g^\top\bar d=0$. Hence $\bar d\in T_\Delta(p)\cap g^\perp=\mathcal K$.

Next, we show that each cluster point solves~\eqref{eq:response-qp}. Fix an arbitrary $z\in\mathcal K$. For indices with $p_i=0$, $z_i\geq0$; for the remaining finitely many indices, $p_i>0$ ensures $p_i+\theta z_i\geq0$ for sufficiently small $\theta>0$. Also $\one^\top z=0$. Thus $p+\theta z\in\Delta$. Moreover, continuity of the Gramian gives
\[
W_{A(\theta)}(p+\theta z)\to W_A(p)\succ0.
\]
Hence, \(p+\theta z\in X_{A(\theta)}\) for sufficiently small \(\theta>0\).

By Taylor's formula for the allocation variable,
\[
 \begin{aligned}
 &h_\theta(p+\theta v)-h_\theta(p)\\
 &=\theta\nabla_u h_\theta(p)^\top v
 +\theta^2\int_0^1(1-t)v^\top
 \nabla^2_{uu}h_\theta(p+t\theta v)v\,dt.
 \end{aligned}
\]
For $v$ in any fixed bounded set, all points in this expression are in a common positive-definite neighborhood for small $\theta$. Moreover, by the definition of $r$, $\nabla_u h_\theta(p)=g+\theta r+o(\theta)$, and the continuity of the Hessian gives $\nabla^2_{uu}h_\theta(p+t\theta v)=H+o(1)$ uniformly for $t\in[0,1]$ and $v$ in bounded sets. Consequently,
\begin{equation}\label{eq:taylor}
 h_\theta(p+\theta v)-h_\theta(p)
 =\theta g^\top v+\theta^2\phi(v)+o(\theta^2),
\end{equation}
\begin{equation}\label{eq:phi}
 \phi(v)=\tfrac12v^\top Hv+r^\top v,
\end{equation}
with a remainder uniform on that bounded set.

Apply \eqref{eq:taylor} to the bounded family $d_\theta$ and to the fixed $z$. By optimality of \(p_\theta\), we have \(h_\theta(p_\theta)\le h_\theta(p+\theta z)\). Applying~\eqref{eq:taylor} with $v=d_\theta$ and $v=z$,
and using $p_\theta=p+\theta d_\theta$
and $g^\top z=0$, we obtain
\begin{equation}\label{eq:optimality-comparison}
 \theta g^\top d_\theta+\theta^2\phi(d_\theta)
 \leq\theta^2\phi(z)+o(\theta^2).
\end{equation}
Since $g^\top d_\theta\geq0$, dividing \eqref{eq:optimality-comparison} by $\theta^2$ gives $\phi(d_\theta)\leq\phi(z)+o(1)$. Along the subsequence that converges to $\bar d$, the continuity of $\phi$ gives $\phi(\bar d)\leq\phi(z)$. Since the comparison holds for every fixed point $z\in\mathcal K$, and $\bar d\in\mathcal K$, every cluster point solves \eqref{eq:response-qp}.

Finally, $\mathcal K$ is a nonempty closed convex cone contained in
$\one^\perp$, and $\phi$ is strongly convex on $\one^\perp$
by Lemma~\ref{lem:curvature}. Hence \eqref{eq:response-qp} has a unique minimizer $d$. Since the family
$d_\theta$ is bounded and each of its cluster points solves \eqref{eq:response-qp},
$d_\theta \to d$, which proves \eqref{eq:score-expansion}.
\hfill$\square$

Under strict complementarity, the critical-cone QP~\eqref{eq:response-qp} reduces to a linear system. Recall the index sets $J$ and $I_0$ defined in~\eqref{eq:index_set}, the Hessian
$H=H_{\rm CS}(p,A)$, and the gradient derivative
$r$ defined in~\eqref{eq:forcing-definition}.
When $I_0=\emptyset$, the critical cone~\eqref{eq:critical-cone}
contains only equality constraints, leading to
the following result.

\begin{corollary}[Fixed-face reduction]\label{cor:fixed-face}
If strict complementarity holds, so that $I_0=\varnothing$, then $d_{J^c}=0$, and $d_J$ is the unique solution of
\begin{equation}\label{eq:fixed-face}
 \begin{bmatrix}H_{JJ}&\one\\\one^\top&0\end{bmatrix}
 \begin{bmatrix}d_J\\\lambda\end{bmatrix}
 =\begin{bmatrix}-r_J\\0\end{bmatrix}.
\end{equation}
\end{corollary}
\begin{IEEEproof}
Strict complementarity means that every zero component of $p$ has a positive multiplier, so $J^c=I_+$ and $I_0=\varnothing$. Equation~\eqref{eq:critical-cone} therefore becomes the linear subspace $\{z:z_{J^c}=0,\one^\top z_J=0\}$. Substituting $z_{J^c}=0$ into the QP leaves the equality-constrained objective $\tfrac12z_J^\top H_{JJ}z_J+r_J^\top z_J$. Its stationarity and feasibility equations are $H_{JJ}d_J+r_J+\lambda\one=0$ and $\one^\top d_J=0$, which are precisely \eqref{eq:fixed-face}. Convexity makes these conditions sufficient as well as necessary.

To prove nonsingularity, let $(u,\alpha)$ belong to the kernel of the matrix in \eqref{eq:fixed-face}. Then $H_{JJ}u+\alpha\one=0$ and $\one^\top u=0$. Multiplying the first equation by $u^\top$ gives $u^\top H_{JJ}u=0$. Extend $u$ by zeros on $J^c$; the resulting vector lies in $\one^\perp$, so Lemma~\ref{lem:curvature} forces $u=0$. Since $J$ is nonempty ($p\in\Delta$), the remaining equation $\alpha\one=0$ gives $\alpha=0$. The kernel is trivial, hence the square saddle-point matrix is nonsingular and the solution is unique.
\end{IEEEproof}

\section{Analytical comparison of VCS and AECS local responses}
\label{app:comparison}
Although $f_{\mathrm V}$ and $f_{\mathrm A}$ depend on the same
node-wise Gramians~\eqref{eq:Wi}, their curvatures \eqref{eq:HV}
and~\eqref{eq:HA}, gradient derivatives \eqref{eq:forcing-V}
and~\eqref{eq:forcing-A}, and the supports of their optimal
allocations, which determine the critical cone~\eqref{eq:critical-cone},
can differ. The following two-node example shows that neither VCS nor AECS has uniformly larger score response norms under fractional edge weakening.

\begin{proposition}[No universal ordering of VCS and AECS response norms]\label{prop:ordering}
For $a,b,c>0$, consider the Metzler Hurwitz matrix
\begin{equation}\label{eq:Ab}
 A_b=\begin{bmatrix}-a&b\\0&-c\end{bmatrix}.
\end{equation}
Along weakening of its edge $2\to1$, for $\rm CS \in \{V, A\}$, define 
\[
d^{\rm CS}(b)=
\left.
\frac{d}{d\theta^+}
p^{\rm CS}(A_b-\theta bE_{12})
\right|_{\theta=0}.
\]
Then,
\begin{equation}\label{eq:S-V-small}
\|d^{\rm V}(b)\|_2=\frac{\sqrt2\,b^2}{(a+c)^2}+O(b^4),
\end{equation}
\begin{equation}\label{eq:S-A-small}
 \|d^{\rm A}(b)\|_2
 =\frac{2\sqrt2\,cb^2}{(\sqrt a+\sqrt c)^2(a+c)^2}+O(b^4).
\end{equation}
Consequently,
\begin{equation}\label{eq:ratio}
 \lim_{b\downarrow0}\frac{\|d^{\rm A}(b)\|_2}{\|d^{\rm V}(b)\|_2}
 =\frac{2c}{(\sqrt a+\sqrt c)^2}.
\end{equation}
Neither $\|d^{\rm V}(b)\|_2\geq \|d^{\rm A}(b)\|_2$ nor its reverse holds for every Metzler Hurwitz system and edge.
\end{proposition}

\begin{IEEEproof}
Solving \eqref{eq:lyapunov} for $A_b$ gives
\begin{equation}\label{eq:example-Gramians}
 \begin{aligned}
 W_1&=\begin{bmatrix}(2a)^{-1}&0\\0&0\end{bmatrix},\\
 W_2&=\begin{bmatrix}
 \dfrac{b^2}{2ac(a+c)}&\dfrac{b}{2c(a+c)}\\[1ex]
 \dfrac{b}{2c(a+c)}&\dfrac1{2c}
 \end{bmatrix}.
 \end{aligned}
\end{equation}
For $p=(x,1-x)^\top$ and $0<x<1$, we have
\begin{equation}\label{eq:example-determinant}
 \det W_{A_b}(p)
 =\frac{(1-x)[b^2+((a+c)^2-b^2)x]}{4ac(a+c)^2}.
\end{equation}

Let $\sigma=a+c$. Maximizing \eqref{eq:example-determinant}, a strictly concave quadratic in $x$ for $b^2<\sigma^2/2$, gives
\begin{equation}\label{eq:xV}
 x_{\mathrm V}(b)=\frac{(a+c)^2-2b^2}{2((a+c)^2-b^2)}=\frac12-\frac{b^2}{2(a+c)^2}+O(b^4).
\end{equation}
Differentiating the first equality in~\eqref{eq:xV} gives
\[
 x_{\mathrm V}'(b)=-\frac{b\sigma^2}{(\sigma^2-b^2)^2}
 =-\frac{b}{\sigma^2}+O(b^3).
\]

A direct computation yields the AECS objective
\begin{equation}\label{eq:example-A-objective}
 f_{\mathrm A}(x,A_b)
 =\frac{2(a+c)[a^2(1-x)+ac+c^2x+b^2(1-x)]}
 {(1-x)[(a+c)^2x+b^2(1-x)]}.
\end{equation}
At $b=0$, this objective reduces to $2a/x+2c/(1-x)$, whose unique minimizer is $x_0=\sqrt a/(\sqrt a+\sqrt c)$.

Set $t=b^2$ and denote the rational objective in \eqref{eq:example-A-objective} by $F(x,t)$. Its denominator is nonzero in a neighborhood of $(x_0,0)$, so $F$ is smooth there. Because $F_x(x_0,0)=0$ and $F_{xx}(x_0,0)>0$, the implicit function theorem applied to $F_x=0$ gives a unique smooth stationary branch $\xi(t)$ with $\xi(0)=x_0$. It remains in $(0,1)$ for small $t$. For $t\geq0$, Lemma~\ref{lem:curvature}, applied to the allocation direction $(1,-1)^\top$, makes the objective strictly convex on its feasible interval. Hence $\xi(t)$ is the unique global minimizer. Thus $x_{\mathrm A}(b)=\xi(b^2)$.

Differentiating~\eqref{eq:example-A-objective} with respect to $t$ at $t=0$ gives
\[
F_{xt}(x,0)=\frac{4a(1-x)}{(a+c)^2x^3}.
\]
Evaluating at $x_0$ gives
\begin{equation}\label{eq:dxAdt}
 \xi'(0)
 =-\frac{F_{xt}(x_0,0)}{F_{xx}(x_0,0)}
 =-\frac{c}{(\sqrt a+\sqrt c)^2(a+c)^2}.
\end{equation}
Smoothness of $\xi$ gives $\xi'(t)=\xi'(0)+O(t)$. Applying the chain rule to $x_{\mathrm A}(b)=\xi(b^2)$ therefore gives
\begin{equation}\label{eq:xA-derivative}
 x_{\mathrm A}'(b)
 =-\frac{2cb}{(\sqrt a+\sqrt c)^2(a+c)^2}+O(b^3).
\end{equation}
Under fractional weakening of the edge $2\to1$,
the weight $b$ changes to $b(1-\theta)$.
Differentiating the optimal allocation with respect
to $\theta$ at $\theta=0$ gives
\[
d^{\rm CS}(b)
=
\begin{pmatrix}
-bx_{\rm CS}'(b)\\
 bx_{\rm CS}'(b)
\end{pmatrix},
\qquad
\|d^{\rm CS}(b)\|_2
=
\sqrt{2}\,b|x_{\rm CS}'(b)|.
\]
Substituting the expansions of $x_{\rm V}'(b)$
and $x_{\rm A}'(b)$ yields~\eqref{eq:S-V-small} and~\eqref{eq:S-A-small},
respectively.

Since the leading coefficient in~\eqref{eq:S-V-small} is positive,
dividing~\eqref{eq:S-A-small} by~\eqref{eq:S-V-small} gives~\eqref{eq:ratio}.
The limiting ratio equals $1/2$ when $a=c$,
whereas it exceeds one when
$c/a>3+2\sqrt{2}$.
These two choices give opposite strict orderings
of the response norms for sufficiently small $b>0$,
proving Proposition~\ref{prop:ordering}.
\end{IEEEproof}

\section{Admissibility along a full Metzler deletion path}
\label{app:metzler}

A real matrix is Metzler if its off-diagonal entries are nonnegative. For this class, reducing a single off-diagonal entry preserves the stability required by the infinite-horizon formulation.

\begin{theorem}[Metzler deletion paths]\label{thm:metzler}
Let $A$ be Metzler and Hurwitz, with $A_{rs}>0$, $r\ne s$. Every matrix $A(\theta)$ in \eqref{eq:deletion} is Metzler and Hurwitz for $0\leq\theta\leq1$.
\end{theorem}

\begin{IEEEproof}
For $0\leq\theta\leq1$, the only changed off-diagonal entry is $(1-\theta)A_{rs}\geq0$, so $A(\theta)$ remains Metzler.
By \cite[Prop.~1]{rantzer2015scalable}, a Metzler matrix $B$ is Hurwitz if and only if there exists a vector $v>0$ such that $Bv<0$, where the inequalities are entrywise.
Since $A$ is Metzler and Hurwitz, choose $v>0$ with $Av<0$. Then
$A(\theta)v
 =Av-\theta A_{rs}v_s e_r
 \leq Av<0$.
The same characterization therefore shows that $A(\theta)$ is Hurwitz for every $0\leq\theta\leq1$.
\end{IEEEproof}

\bibliographystyle{IEEEtran}
\bibliography{reference}

\end{document}